\documentstyle[12pt]{article}

\newcommand{\vs}{\vspace{.2in}}

\newcommand{\noin}{\noindent}

\def\C{{\bf C}}
\def\P{{\bf P}}
\def\R{{\bf R}}
\def\Z{{\bf Z}}

\def\tp{\tilde {p}}

\def\cp1{{\bf CP}^1}
\def\k{\bar {k}}

\def\d{{\rm dim}_{\C}\,}

\def\al3{(\alpha_3 -1 )}

\def\\{\vskip 0.4cm}
\def\is{\simeq}
\def\qed{\hfill $\Box$}
\def\p{\tilde p}
\begin{document}

\centerline {\LARGE On exotic algebraic structures on affine spaces} \\

\vs

\centerline {\bf M. Zaidenberg} \\

\vs

We review the recent developement on the subject, emphasizing its analytic
aspects and pointing out some open problems. The study of exotic $\C^n$-s is at
the very beginning, and hopefully this survey would be useful to learn more
about these unusual and beautiful objects.

All algebraic varieties considered below are usually assumed being smooth,
reduced, irreducible, and defined over $\C$. {\it Isomorphism} means biregular
isomorphism of algebraic varieties and is denoted by $\simeq$. \\

It is my pleasure to thank Sh. Kaliman and P. Russell. Without their friendly
help and advice this survey would not be written. \\
\vs

\centerline {\LARGE Contents}\\

\noin {\bf 1. Product structures}\\

\noin {\bf 2. Kaliman's modification} \\

\noin {\bf 3. Hyperbolic modification} \\

\noin {\bf 4. Dimca's and Kaliman's examples of exotic hypersurfaces} \\

\noin {\bf 5. Russell's $\C^*$--threefolds and the Makar--Limanov invariant} \\

\noin {\bf 6. APPENDIX: Simply connectedness of $\C^*$--equivariant cyclic
coverings}\\

\noin {\bf 7. Concluding remarks} \\

\vs

\centerline {{\bf 1. Product structures}} \\

\noin {\bf 1.1. Definition.} Let $X$ be a smooth affine algebraic variety. It
is called {\it an exotic $\C^n$} if $X$ is diffeomorphic to $\R^{2n}$ but not
isomorphic to $\C^n$.

The following characterization of $\C^2$, due to C. P. Ramanujam [Ram] shows
that there is no exotic $\C^2$. \\

\noin {\bf 1.2. Ramanujam's Theorem.} {\it A smooth affine algebraic surface
$X$ is isomorphic to $\C^2$ iff the groups ${\rm H}_2 (X; \Z),\,\pi_1 (X)$, and
$\pi_1^{\infty} (X)$ all are trivial. In particular, if $X$ is homeomorphic to
$\C^2$, then $X$ is isomorphic to $\C^2$. } \\

\noin Here $\pi_1^{\infty} (X)$ denotes the fundamental group at infinity of
$X$. If $X$ is the interior of a compact real manifiold with the boundary
$\partial X$, then $\pi_1^{\infty} (X) = \pi_1 (\partial X)$.

To recognize an exotic $\C^n$ one has to verify two properties from Definition
1.1. As for the first one, the next criterion could be useful. \\

\noin {\bf 1.3. Ramanujam--Dimca's Theorem} ([Ram; Di], see also [Ka 1, Lemma
1; tD 1, Theorem 3.14]). {\it  Let $X$ be a smooth affine algebraic variety of
complex dimension $n \ge 3$. Then $X$ is diffeomorphic to $\R^{2n}$ iff $X$ is
contractible, or, equivalently, iff $\pi_1 (X) = {\bf 1}$ and ${\rm H}_i (X;
\Z) = 0$ for $i = 1,\dots,n$. } \\

\noin The proof is based on the h-cobordism theorem. The main point is to show
that being contractible, $X$ possesses a smooth simply connected boundary
$\partial X$. The latter follows from the Lefschetz Hyperplane Section Theorem.
In the case when $X$ as below is a product of two contractible varieties,
instead of Lefschetz--type arguments one can apply the Van Kampen Theorem. \\

\noin As follows from [Ram] the statement of Theorem 1.3 does not hold for $n =
2$. Indeed, Ramanujam constructed an example of a smooth contractible affine
surface $S_0$ which is not homeomorphic to $\R^4$ and, moreover, has an
infinite group $\pi_1^{\infty} (S_0) $. See e.g. [GuMi, FlZa, PtD 3] for
further information on such surfaces. \\

\noin  {\bf 1.4.} We recall {\bf Zariski's Cancellation Problem}: \\

\noin {\it Given an isomorphism $X \times \C^k \simeq \C^{n + k}$, does it
follows that $X \simeq \C^n$?} \\

\noin For $n \ge 3$ the problem is still open. For n = 2 the positive answer
was obtained by the Miyanishi-Sugie and Fujita Theorem [MiSu, Fu 1].  This
theorem provides us, in certain cases, with a tool to distinguish exotic
$\C^n$--s. Despite of the fact that it was proven later on, the first examples
of exotic $\C^n$--s for any $n \ge 3$ had been alluded to already in [Ram]:
take $X^n = S_0 \times \C^{n - 2}$, where $S_0$ is the Ramanujam surface. See
[Za 1,2,4] for some further examples of exotic structures of product type. \\

\noin {\bf 1.5. Exotic $\C^n$--s of log-general type.} More generally, let
$X_i, \,i=1,\dots,m$, be smooth contractible affine varieties. Then by
Ramanujam--Dimca's Theorem 1.3 $X := (\prod_{i = 1}^m X_i ) \times \C^r$ is
diffeomorphic to $\R^{2n}$ as soon as $n :={\rm dim}_{\C}\,X $ is at least $3$.
To distinguish such a product structure $X$ from the standard $\C^n$, the
logarithmic Kodaira dimension $\k$ could be available (see [Ii 1]). But, if $r
> 0$, then $\k (X) = \k (\C^n) = - \infty$. So, assume for the moment that $r =
0$. Then $\k (X) = \sum_{i=1}^m \k(X_i) \ge 0$ as far as $\k(X_i) \ge 0$ for
each $i = 1,\dots,m$. Therefore, in this case $X$ is an exotic $\C^n$.

For instance, put $X = (S_0)^m ,\,m > 1$, where $S_0$ is the Ramanujam surface.
Then $X$ is an exotic $\C^n$ of log-general type, i.e. $\k(X) = n$,  where $n =
2m ={\rm dim}_{\C}\,X $ (indeed, by [Ii 2] we have $\k(S_0) = 2$). Using other
contractible surfaces of log-general type as factors, one may construct more
such examples \footnote{this remark is due to T. tom Dieck [tD 1]}, and up to
now these are the only known ones. Note that all of them are of even
dimensions. \\

\noin {\bf 1.6.} {\bf Remark.} Let $S$ be a contractible, or at least acyclic
smooth surface of non--negative Kodaira dimension. Then $\k(S) = 1$ or $2$ by
the Fujita classification of open surfaces with $\k = 0$ [Fu 2]\footnote{by the
Miyanishi-Sugie and Fujita Theorem [MiSu, Fu 1] (see also [Ru 2]) $\C^2$ is the
only contractible, or even acyclic, surface of $\k = -\infty$.}. Thus, for the
product structure $X = \prod_{i = 1}^m S_i$ obtained by means of surfaces $S_i$
as above its log-Kodaira dimension takes the values in the interval $m \le
\k(X) \le 2m = {\rm dim}_{\C}\,X $. Later on we will see that there exist
exotic $\C^n$--s $X$ with $1 \le \k(X) \le [n/2]$. \\

\noin {\bf 1.7. Examples.} A source of examples provide contractible affine
surfaces of log-Kodaira dimension one. The complete list of them was obtained
by R. V. Gurjar and M. Miyanishi [GuMi]. T. tom Dieck and T. Petrie [PtD 1]
realized some of them as hypersurfaces in $\C^3$. Namely, put $p_{k,\,l}
(x,\,y,\,z) = ((xz + 1)^k - (yz + 1)^l)/z \in \C [x,\,y,\,z]$, where $(k, l) =
1,\, k, \,l \ge 2$. Then all the fibres of the polynomial $p_{k,\,l}$ are
contractible surfaces in $\C^3$. All its non-zero fibres are smooth surfaces of
log-Kodaira dimension one, isomorphic to $S_{k,\,l} = p_{k,\,l}^{-1} (1)$ (the
zero fibre has non--isolated singularities).

Due to Theorem 1.8 below, the product $S_{k,\,l} \times \C^{n - 2}$ is a
hypersurface in $\C^{n + 1}$ which is an exotic $\C^n$. Moreover, all the
fibres of $p_{k,\,l}$ regarded as polynomial on $\C^{n + 1}$ are contractible
hypersurfaces, and all but the zero one are exotic $\C^n$--s.

Recently Sh. Kaliman and L. Makar--Limanov [KaML 3] have shown that all the
log-Kodaira dimension one contractible surfaces from the Gurjar--Miyanishi list
admit embeddings into $\C^3$, thus providing deformation families of
contractible hypersurfaces (see (1.9) below). \\

\noin {\bf Proposition (Sh. Kaliman, L. Makar--Limanov).} {\it Any smooth
contractible surface $S$ of log-Kodaira dimension one can be given as a surface
in $\C^3$ with the equation $p_{k,\,l,\,m,\,f} (x,\,y,\,z) = 0$, where
$$p_{k,\,l,\,m,\,f} (x,\,y,\,z) = [(z^m x + f(z))^k - (z^m y + g(z))^l -
z]/z^m\,\,,$$ $(k,\,l) = 1, \,k,\,l \ge 2, \,m \ge 1,\,f,\,g \in \C[z],
\,\,{\rm deg}\,f,\,{\rm deg}\,g < m,\,f(0) = g(0) = 1$ and where $g$ is
uniquely defined by $f$ (which can be taken arbitrary) in view of the
assumption that $p_{k,\,l,\,m,\,f}$ is a polynomial. }\\

\noin {\bf 1.8.} Besides the log-Kodaira dimension, one may equally use other
invariants to distinguish exotic product--structures, for instance the
logarithmic plurigenera ${\bar P}_{m_1,\dots,m_l}$ [Ii 1,2]. Indeed, if $X :=
(\prod_{i = 1}^m X_i )$ and at least one of the factors $X_i$ has a non--zero
log--plurigenus, then by K\"unneth's Formula the corresponding log--plurigenus
of $X$ does not vanish. Moreover, we have the following theorem. \\

\noin {\bf Iitaka--Fujita's Strong Cancellation Theorem} [IiFu]. {\it Let
$X,\,Y,\,A_1,\,A_2$ be algebraic varieties such that ${\rm dim}_{\C}\,X = {\rm
dim}_{\C}\,Y,\,\k(X) \ge 0$, and the log--plurigenera ${\bar P}_{m_1,\dots,m_l}
(A_i)$ all vanish, $i = 1,2$. Then, given an isomorphism $\Phi \,:\,X \times
A_1 \to Y \times A_2$, there exists an isomorphism $\varphi\,:\,X \to Y$ making
the following diagram commutative: \\

\begin{picture}(200,80)
\unitlength0.2em
\thicklines
\put(62,25){$X \times A_1$}
\put(105,25){$Y \times A_2$}
\put(84,27){$\vector(1,0){15}$}
\put(66,5){$X$}
\put(84,6){$\vector(1,0){15}$}
\put(109,5){$Y$}
\put(69,22){$\vector(0,-1){11}$}
\put(111,22){$\vector(0,-1){11}$}
\put(91,10){$\varphi$}
\put(91,30){$\Phi$}
\put(160,11){$(1)$}
\end{picture}\\
\noin (the vertical arrows are the canonical projections). } \\

\noin Applying this theorem to the case when $A_1 = A_2 = \C^r$ we see that, as
soon as $\k(X) \ge 0$, the factor $\C^r$ cancells. Therefore, if $X$ and $Y$
are contractible varieties such that $n :=\d X + r \ge 3$ and $0 \le \k(X) \neq
\k(Y)$, then  $X \times \C^r$ is an exotic $\C^n$ non--isomorphic to $Y \times
\C^r$. This allows us to distinguish some exotic product structures of negative
log--Kodaira dimension.  \\

\noin {\bf 1.9. Deformable exotic product structures} \\

\noin {\bf Theorem (Flenner--Zaidenberg)} [FlZa]. {\it For any $n \ge 3$ and $m
\ge 1$ there exists a family of versal deformations ${\tilde f}\,:\,{\tilde
X}^{n + m} \to B^m$ of exotic product structures on $\C^n$ over a smooth
quasiprojective base $B^m$ of dimension $m$.} \\

\noin Here $\tilde f$ is a smooth morphism such that for any $b \in B^m$ the
fibre $X_b = {\tilde f}^{-1} (b)$ is an exotic product structure on $\C^n$ and
$X_b \not\simeq X_{b'}$ if $b \neq b'$ are generic points of $B^m$.

The construction uses families of versal deformations of contractible surfaces
of log--Kodaira dimension one. If $ f\,:\,{\tilde S}^{2 + m} \to B^m$ is such a
family, then by the Iitaka--Fujita Theorem 1.8 we can take ${\tilde X}^{n + m}
= {\tilde S}^{2 + m} \times \C^{n - 2}$. Due to [KaML 3] (see 1.7) this leads
to families of exotic hypersurfaces in $\C^{n+1}$.

In the next section we will give another construction of deformable exotic
$\C^n$--s due to Sh. Kaliman. But still we do not know the answer to the
following \\

\noin {\bf Question.} {\it Does there exist a deformable exotic $\C^n$ of
log--general type?}\footnote{cf. the Rigidity Conjecture for acyclic surfaces
of  log--general type in [FlZa].}\\

\noin {\bf 1.10. Analytically exotic $\C^n$--s.} We say that an exotic $\C^n$
is {\it analytically exotic} if it is not even biholomorphic to $\C^n$. The
following result shows that among the product structures constructed above
there are    analytically exotic ones. \\

\noin {\bf Strong Analytic Cancellation Theorem} [Za 2]. {\it Let $X,\,Y$ be
smooth quasi--projective varieties of log--general type. If $\Phi\,:\, X \times
\C^r \to Y \times \C^m$ is a biholomorphism, then $m = r$ and there exists an
isomorphism $\varphi\,:\,X \to Y$ which makes diagram $(1)$ commutative, where
$A_i, \, i=1,2$ are replaced by $\C^r$. } \\

\noin Thus, if $S$ is a contractible surface of log--general type and $n \ge
3$, then $X = S \times \C^{n - 2}$ is an analytically exotic $\C^n$. The same
is true for the product structures $X = S \times M^{n - 2}$, where $S$ is as
above and $M^{n - 2}$ is any contractible affine variety of dimension $n - 2$.
Indeed,  for such an $X$ its Eisenman--Kobayashi intrinsic 2--measure form
$E^{(2)}_X$ does not vanish identically (this useful remark is due to Sh.
Kaliman [Ka 2]). More generally, one can consider the maximal value of $k$ for
which $E^{(k)}_X$ does not vanish identically. This yields a coarse analytic
invariant which replaces the log--Kodaira dimension as it has been used in 1.8
above, and so it permits to distinguish certain analytically exotic $\C^n$-s up
to a biholomorphism.

Another remarkable property of the above exotic product structures on $\C^n$ is
that they contain no copy of $\C^{n-1}$. \\

\noin {\bf 1.11. Theorem} [Za 2; Ka 2]. {\it Let $X = S \times \C^{n - 2}$ be
an exotic product structure on $\C^n$, where $S$ is a contractible surface of
log-general type. Then there is no regular injection $\C^{n - 1} \to X$; in
particular, there is no algebraic hypersurface in $X$ isomorphic to $\C^{n -
1}$. Moreover, there is no proper holomorphic injection $\C^{n - 1} \to X$, and
therefore there is no closed analytic hypersurface in $X$ biholomorphic to
$\C^{n - 1}$. } \\

Theorem 1.10 shows that the following is likely to be true. \\

\noin {\bf 1.12. Conjecture.} {\it Any exotic $\C^n$ is analytically exotic.}
\\

\noin {\bf 1.13. Problem.} {\it Does there exist a pair of biholomorphic but
not isomorphic exotic $\C^n$-s? Does there exist a non-trivial deformation
family of exotic $\C^n$-s with the same underlying analytic structure?} \\

\noin We even do not know whether the deformation families of exotic
product--structures constructed in the proof of Theorem 1.9 are versal in the
analytic sense. To this point, the knowledge of the collection of entire curves
(i.e. holomorphic images of $\C$) in the contractible surfaces of log--Kodaira
dimension one could be useful. What is the set of tangent directions of such
curves in the tangent bundle? Does the degeneration locus of the Kobayashi
pseudo--distance provide a non-trivial analytic invariant of such surfaces? \\

\centerline {\bf 2. Kaliman's modification} \\

\noin {\bf 2.1. Definition} (cf. [Ka 2]). Consider a triple $(X,\,H,\,C)$
consisting of an algebraic variety $X$, an irreducible hypersurface $H$ in $X$
and a closed subvariety $C$ of $H$ with ${\rm codim}_X C \ge 2$. Let
$\sigma_C\,:\,{\tilde X} \to X$ be the blow-up of $X$ at the ideal sheaf of
$C$. Let $E \subset {\tilde X}$ be the exceptional divisor of $\sigma_C$ and
${\tilde H} \subset {\tilde X}$ be the proper transform of $H$. {\it The
Kaliman modification}\footnote{in the surface case a similar transform was
called a half--point attachement (detachment) in [Fu 2].} consists in replacing
the triple  $(X,\,H,\,C)$ by the pair $(X',\,E')$, where $X' =  {\tilde X}
\setminus {\tilde H}$ and  $E' =  E \setminus {\tilde H}$. We will also say
that {\it $X'$ is the Kaliman modification of $X$ along $H$ with center $C$}.

A triple $(X,\,H,\,C)$ resp. a pair $(X',\,E')$ as above will be called {\it a
smooth contractible affine triple} resp. {\it a smooth contractible affine
pair} if all its members are smooth contractible affine varieties. \\

\noin {\bf 2.2. Theorem} [Ka 2, Theorem 3.5]. {\it The Kaliman modification of
a smooth contractible affine triple is a smooth contractible affine pair.} \\

\noin The statement of the theorem remains valid under the assumption that the
hypersurface $H$ (not necessarily smooth any more) is a toplogical cell and $C
\subset {\rm reg}\,H$, while all other conditions on $X$ and $C$ being
preserved (see  [Ka 2, Theorem 3.5]). However, we do not know whether the
smoothness of $C$ is essential\footnote{See [Ka 2] for an example which shows
that the theorem does not work without the condition $C \subset {\rm reg}\,H$,
even with a smooth $C$.}. \\

\noin {\bf 2.3. Examples.}  The Kaliman modification produces new examples of
analytically exotic $\C^n$--s and of their versal deformations. Let
$(X,\,H,\,C)$ be a smooth contractible affine triple, where $X$ is an exotic
$\C^n$ such that certain Eisenman--Kobayashi form $E^{(k)}_X$ is different from
zero at the points of an open subset $U \subset X$. Performing the Kaliman
modification we arrive again at an exotic $\C^n$, call it $X'$, which has a
non-zero form $E^{(k)}_{X'}$. Indeed, by Theorems 1.3 and 2.2 $X'$ is
diffeomorphic to $\R^{2n}$. Furthermore, the restriction $
\sigma_C\,|\,X'\,:\,X' \to X$ is a dominant holomorphic mapping which is a
contraction with respect to the Eisenman--Kobayashi forms, i.e. $\sigma_C^* \,
E^{(k)}_X \le E^{(k)}_{X'}$, and whence $E^{(k)}_{X'} \not\equiv 0$. Thus, $X'$
is an analytically exotic $\C^n$.

For instance, fix a point $s_0$ in a smooth contractible affine surface $S$ of
log--general type and put $X = S \times \C^{n-2} ,\,H = S \times \C^{n-3}$ and
$C = \{s_0\} \times \C^k$, where $0 \le k \le n-3$. By Sakai's Theorem [Sak]
$S$ is measure hyperbolic, i.e.  $E^{(2)}_S$ is positive on a subset of $S$
whose complement has measure zero. Therefore, $E^{(2)}_X = {\rm pr}^*_S\,
E^{(2)}_S$ is different from zero at the points of a massive subset of $X$
(where ${\rm pr}_S\,:\,X \to S$ is the natural projection). Performing the
Kaliman modification, by Theorem 2.2, we obtain a smooth contractible pair
$(X',\,E')$, where $X'$ is an exotic $\C^n$ which have a non-zero form
$E^{(2)}_{X'}$. It is easily seen that $E' \simeq \C^{n-1}$. Hence, $X'$ is not
biholomorphic to any ${\tilde X} = {\tilde S} \times \C^{n-2}$ as above.
Indeed, by Theorem 1.11 such an ${\tilde X}$ does not contain any biholomorphic
image of $\C^{n-1}$. And also it is not  biholomorphic to any exotic ${\tilde
X} = {\tilde S} \times \C^{n-2
}$, where $\k ({\tilde S}) =  1$, because for the latter product the form
$E^{(2)}_{\tilde X}$ vanishes at a Zariski open subset. \footnote{we are
grateful to Sh. Kaliman for this remark.}
This proves the following \\

\noin {\bf 2.4. Proposition.} {\it For $n = 3\,$ $X'$ as above is an exotic
$\C^3\,$ which is not biholomorphic to any exotic product--structure on
$\C^3$.} \\

\noin {\bf 2.5. Deformable analytically exotic $\C^n$} \\

\noin By iterating the construction used in the preceeding example, Sh. Kaliman
 obtained the following result (cf. Theorem 1.9 above). \\

\noin {\bf Theorem} [Ka 2, sect.4]. {\it For any $n \ge 3$ there exist versal
deformation families of analytically  exotic $\C^n$--s with any given number of
moduli}. \\

\noin The proof proceeds as follows. Start with an exotic $\C^n\,$, $X = S
\times \C^{n-2}$, as above. Fix $m$ distinct points $s_1,\dots, s_m \in S$ and
$m$ disjoint affine hyperplanes $A_1,\dots, A_m$ in $\C^{n-2}$. Put $H_i = S
\times A_i,\,C = \{s_i\} \times A_i$, and perform the Kaliman modifications
along $H_i$ with centers $C_i$ for $i=1,\dots,m$. Then we result with a family
$X' = X'\,(s_1,\dots, s_m ,\,A_1,\dots, A_m)$ of analytically  exotic $\C^n$--s
endowed each one with $m$ disjoint hypersurfaces $E'_1,\dots,E'_m$ isomorphic
to $\C^{n-1}$. And they are the only biholomorphic images of $\C^{n-1}$ in $X'$
for fixed data [Ka 2, Lemma 4.1]. Now it is not difficult to check that the
positions of the points $s_1,\dots, s_m$ in $S$ and of the hyperplanes
$A_1,\dots, A_m$ in $\C^{n-2}$ provide the moduli of these exotic structures
considered up to a biholomorphism. \\

\newpage

\centerline {\bf EXOTIC AFFINE HYPERSURFACES} \\

In Sections 3--5 below we review some explicit constructions of contractible
hypersurfaces in $\C^{n+1}$. We discuss different approaches to the recognition
problem for exotic hypersurfaces.

By the Abhyankar--Moh and Suzuki Theorem [AM, Su] the only smooth irreducible
simply connected curves in $\C^2$ are those obtained from the affine line $x =
0$ by means of polynomial coordinate changes. By the Lin--Zaidenberg Theorem
[LiZa] the only simply connected affine plane curves, up to the action of the
group of biregular automorphisms ${\rm Aut}\,\C^2$, are the quasihomogeneous
ones. In particular, each irreducible singular such curve is equivalent to one
and the only one from the sequence $\Gamma_{k,\,l} = \{x^k - y^l =
0\},\,(k,\,l) = 1, \,k>l\ge 2$. Thus, in the case $n = 2$ this classifies
completely the contractible hypersurfaces in $\C^n$.

In contrast, we have already seen in Section $1$ above that there are even
deformation families of smooth contractible surfaces of log--Kodaira dimension
one in $\C^3$, and so for any $n \ge 3$ there are deformation families of
hypersurfaces in $\C^{n+1}$ which are exotic $\C^n$--s. They are far from being
classified in any sense. In particular, {\it no exotic hypersurface in
$\C^{n+1}$ of log--general type is known.} \\

\centerline {\bf 3. Hyperbolic modification} \\

In [tD 1] T. tom Dieck introduced a general construction which, under certain
conditions, represents a given topological resp. complex manifold $Z$ with
possible singularities as an algebraic quotient of an action of a real resp.
complex Lie group $G$ on another such manifold $X$, canonically defined by Z
and a given representation of $G$. This representation should have the unique
fixed point, which should be of hyperbolic type, and so the correspondence $Z
\longmapsto X$ was called {\it the hyperbolic modification}.

When $Z \subset \C^n$ is an affine algebraic variety and $G = \C^*$, the
hyperbolic modification $X$ of $Z$ is an affine algebraic variety in $\C^{n+1}$
endowed with a $\C^*$--action, and $Z \is X//\C^*$. The $n$--th iterate $X_n$
of the hyperbolic modification of $Z$ is endowed with an action of the
$n$--torus $T_n = (\C^* )^n$, and $Z \is X_n//T_n$.

The main advantage of this transform is that it leads, in the case of
hypersurfaces, to new amazing examples of exotic families. \\

\noin {\bf Proposition 3.1} [tD 1, KaML 3]. {\it For any $n \ge 3$ there exist
effectively defined polynomials $p_{k,\,l}^{(n)}$ on  $\C^{n+1}$, where
$(k,\,l) = 1,\,k > l \ge 2$, such that all the fibres $(p_{k,\,l}^{(n)})^{-1}
(c),\,c \in \C$, are exotic $\C^n$--s.} \\

Since we are interested mainly in the hypersurface case, and it is simpler, we
give the precise definition only in this case. \\

\noin {\bf 3.2. Definition} [tD 1]. Let $p \in \C [x_1,\dots,x_n] \setminus
\{0\}$ be a polynomial on $\C^n$ such that $p({\bar 0}) = 0$. The polynomial
$\tp (x, z) = p(xz)/z \in \C[x,\,z]$ on $\C^{n+1}$ is called {\it the
hyperbolic modification of $p$}, and its zero fibre $X_0 = \tp^{-1} (0) \subset
\C^{n+1}$ is called {\it the hyperbolic modification of the zero fibre $Z_0 =
p^{-1} (0)  \subset \C^n$ of $p$}. \\

In the case of the simply connected curves $\Gamma_{k,\,l} = \{x^k - y^l = 0\}$
the hyperbolic modification was already used in [PtD 1]. The origin having been
placed at the smooth point $(1,\,1) \in \Gamma_{k,\,l}$, the hyperbolic
modification gives rise to the Petrie--tom Dieck polynomials $p_{k,\,l} = ((xz
+ 1)^k - (yz + 1)^l)/z$ (see 1.6; see also [PtD 2] for some related
constructions). \\

\noin {\bf 3.3. Some properties of the hyperbolic modification.} Let $\tp (x,
z)$ be the hyperbolic modification of a polynomial $p(x)$ on $\C^n$. Consider
the $\C^*$--action $G(\lambda,\, x, \,z) = (\lambda x,\,\lambda^{-1} z)$ on
$\C^{n+1}$. It is easily seen that $\tp$ is a quasi--invariant of $G$ of weight
$1$, i.e. $\tp \circ G_{\lambda} = \lambda \tp ,\,\lambda \in \C^*$. Denote
$X_1 = \tp^{-1} (1) \subset \C^{n+1}$. Then the restriction $G \,|\, (\C^*
\times X_1)$ yields an isomorphism $\C^* \times X_1 \simeq \C^{n+1} \setminus
X_0$. In particular, $X_1$ is a smooth hypersurface, and all the fibres $X_c =
\tp^{-1} (c),\,c \in \C \setminus \{0\}$, are isomorphic to $X_1$.

Thus, being applied to a hypersurface $Z_0$ in $\C^n$, the hyperbolic
modification  produces actually a pair of distinct hypersurfaces $X_0$ and
$X_1$ in $\C^{n+1}$. The zero fibre $X_0$ inherits the $\C^*$--action
$G\,|\,X_0$. The ring of $G$--invariants coincides with the  subring $\C[zx_1
,\dots, zx_n] \subset \C[x, \,z]$, and $\pi\,:\,\C^{n+1} \ni (x, \,z)
\longmapsto xz \in \C^n$ is the canonical quotient morphism, as well as the
restriction $\pi \,|\,X_0 \,:\,X_0 \to Z_0$. Thus, $Z_0 = X_0 //\C^*$ is the
algebraic quotient. Note that $X_0$ is smooth iff $Z_0$ is so. \\

For the morphism $\pi\,:\,X_0 \to Z_0$ the Kawamata Addition Theorem [Kaw] and
the Iitaka inequality [Ii 1] imply that $\k (Z_0) \le \k (X_0) \le \d Z_0$. The
same holds for any iterated hyperbolic modification $X_0^{(n)}$ of $Z_0$.   \\

The restriction $\pi \,|\,X_1 \,:\,X_1 \to \C^n$ is a birational morphism. From
Proposition 3.6 in [tD 1] it follows that the generic fibre $X_1$ is the
Kaliman modification of $\C^n$ along $Z_0$ with center at the origin (see 2.1).
Combining several statements from [tD 1, (1.1), (1.3), (2.1), (3.1); Ka 2,
(3.5)] we obtain the following \footnote{Another approach to the proof of (a
certain part of) this theorem, based on Kempf--Ness and Neeman results on
algebraic group actions (see e.g. [Kr 1, \S 4]  and references therein) was
proposed in the lectures of T. Petrie at the  Workshop on Open Algebraic
Varieties, CRM, Montr\'eal, December 5--9, 1994; see also [Pe].}\\

\noin {\bf 3.4. Theorem.} {\it Let $Z_0 = p^* (0)$, where $p \in
\C[x_1,\dots,x_n]$. Assume that $Z_0$ is an irreducible reduced divisor on
$\C^n$, which contains the origin and is smooth at the origin. Let $X_0$ resp.
$X_1 \subset \C^{n+1}$ be the zero fibre resp. the generic fibre of the
hyperbolic modification $\tp$ of the polynomial $p$.

\noin a) Let $Z_0$ be smooth. Then both the hyperbolic modification $X_0$ and
the Kaliman modification $X_1$ of $Z_0$ are acyclic resp. contractible as soon
as $Z_0$ is acyclic resp. contractible.

\noin b) Let $Z_0$ be  a topological manifold and has at most isolated
singularities. Then $X_1$ is acyclic resp. contractible if $Z_0$ is acyclic
resp. a topological cell.} \\

This leads to the following result [tD 1, Theorem 3.12].\\

\noin {\bf 3.5. Theorem (T. tom Dieck).} {\it If $Z_0 = p^{-1}(0)$ is a smooth
contractible hypersurface in $\C^n,\, n\ge 3$, then all the fibres of the
hyperbolic modification $\tp$ of $p$ are smooth hypersurfaces in $\C^{n+1}$
diffeomorphic to $\R^{2n}$. If, furthermore, $\k (Z_0 ) \ge 0$, then the zero
fibre $X_0$ of $\tp$ is an exotic $\C^n$ of non--negative log--Kodaira
dimension.} \\

Starting with the Petrie--tom Dieck surface $S_{k,\,l}$ in $\C^3$ of
log--Kodaira dimension one (see 1.7) and iterating the hyperbolic modification,
for any given $n \ge 3$ one can effectively find a polynomial $p^{(n)}_{k,\,l}$
on $\C^{n+1}$ such that all its fibres are smooth and diffeomorphic to
$\R^{2n}$, and the zero fibre is an exotic $\C^n$. In fact, in this particular
case all of them are exotic $\C^n$--s, as has been recently shown in [KaML 3].
More precisely, it was shown that none of these fibres  is dominated by $\C^n$.
This proves Proposition 3.1 above. The Brieskorn--Pham polynomials provide
another examples of this type (in this case one has to apply (b) of Theorem
3.4; see [tD 1, Section 4]).

To manage the general case, it would be useful to prove the following
conjecture, which seems to be interesting by itself. It is easily checked for
$n = 2$. \\

\noin {\bf 3.6. Conjecture. } {\it Let $X_0$ be a special fibre and $X_1$ be a
generic fibre of a primitive polynomial $p \in \C[x_1,\dots,x_n]$. Let $X'_0$
be the desingularization of an irreducible component of $X_0$. Then $\k (X'_0)
\le \k (X_1)$. In other words, the log--Kodaira dimension is lower
semi--continuous on the fibres of a polynomial in $\C^n$. } \\

Note that for $n \ge 5$ the exotic $\C^n$ which is the zero fibre $X_0 =
(p^{(n)}_{k,\,l})^{-1} (0)$ is different from the exotic product structures on
$\C^n$ considered in 1.5 above, since here $1 \le \k (X_0) \le 2$. \\

\begin{center} {\bf 4. Dimca's and Kaliman's examples of exotic hypersurfaces}
\end{center}

\noin As we have seen in the preceeding section, the hyperbolic modification
$\tp$ of a polynomial $p$ on $\C^n$ is a quasi--invariant of weight $1$ of a
linear $\C^*$--action $G$ on $\C^{n+1}$ of mixed type (the latter means that
the $\C^*$--action $G$ has one fixed point only and the linear part of $G$ in
the fixed point has weights of different signs). In the examples considered
below the defining polynomials of exotic hypersurfaces in $\C^{n+1}$ will be
quasi--invariants of weights $> 1$ of regular $\C^*$--actions of mixed type.
The generic fibre $X_1$ of such a polynomial does not need to be contractible
any more. Its zero fibre $X_0$ will be presented as a cyclic branched covering
of $\C^n$ ramified along a hypersurface $Z_0$ with certain properties, which
ensure that $X_0$ is contractible. \\

\noin {\bf 4.1. Dimca's list} [Di, ChoDi]. A. Libgober [Lib] discovered that
the (singular) projective hypersurface $${\bar H}_{n,\,d} = \{x_0^{d-1} x_1 +
x_1^{d-1} x_2 +\dots + x_{n-1}^{d-1} x_n + x_{n+1}^d = 0\} \subset \P^{n+1}$$
has the same homology groups as $\P^n, \,n$ odd. In fact, ${\bar H}_{n,\,d}$ is
a completion of $\C^n$; namely,  $H_{n,\,d} \simeq \C^n$, where $H_{n,\,d} =
{\bar H}_{n,\,d} \setminus \{x_0 = 0\}$.

Generalizing this example, G. Barthel and A. Dimca [BaDi] found some others
homology $\P^n$--s with isolated singularities \footnote{and in particular, all
homology $\P^2$--s which are normal surfaces in $\P^3$ endowed with a
$\C^*$--action.}. These are the projective closures in $\P^{n+1}$ of the affine
hypersurfaces $H_{n,\,d,\,a} \subset \C^{n+1}$ with the equations
$$p_{n,\,d,\,a} (x) = x_1 + x_1^{d-1}x_2 + \dots + x_{n-2}^{d-1} x_{n-1} +
x^{d-a}_{n-1}x_n^a + x^d_{n+1} = 0$$ where $n$ is odd and $1 \le a <
d,\,(a,\,d) = (a,\,d-1) = 1$. \\

\noin {\bf 4.2. Proposition} [Di, Propositions 5, 7; ChoDi, Theorem 5,
Proposition 6]. {\it

\noin a) $H_{n,\,d,\,a} \subset \C^{n+1}$  is diffeomorphic to $\R^{2n}$.

\noin b) For $a = 1$ the map $p_{n,\,d,\,a}\,:\,\C^{n+1} \to \C$ is a smooth
fibre bundle with the fibre diffeomorphic to $\R^{2n}$.

\noin c) For $a > 1$ the generic fibre $X_1 = p_{n,\,d,\,a}^{-1} (1)$ is not
contractible (in fact, its Euler characteristic is different from $1$).

\noin d) $H_{3,\,d,\,1} \simeq \C^3$ and the fibration $p_{3,\,d,\,1}\,:\,\C^4
\to \C$ is algebraically trivial.} \footnote{Sh. Kaliman [Ka 1] noted that for
$d = 3$ it is trivialized by the Nagata automorphism. } \\

A. Dimca posed the question:

\noin {\it Is it true that for $a > 1$ all the hypersurfaces $H_{3,\,d,\,a} =
\{x + x^{d-1}y + y^{d-a}z^a + t^d = 0\} \subset \C^4$ as above are exotic
$\C^3$--s?}

\noin The positive answer was done by [KaML 2], see Theorem 5.10 below. It is
still unknown whether the same is true in higher dimensions. \\

The following criterion of contractibility of cyclic coverings, proposed in [Ka
1, Theorem A] (see also [tD 2]), allows one to establish that the hypersurfaces
like those in the previous examples and more general ones are contractible.\\

\noin {\bf 4.3. Theorem (Kaliman)}. {\it Let a polynomial $q \in \C[x_1,\dots,
x_n]$ be a quasi--invariant of a positive weight $l$ of a regular
$\C^*$--action $G$ on $\C^n$. Suppose that the zero fibre $Z_0 = q^* (0)$ is a
smooth, reduced, and irreducible divisor in $\C^n$ such that

\noin i)  $\pi_1 (\C^n \setminus Z_0 ) \approx \Z$ \footnote{see footnote 11
below.};

\noin ii) for some prime $k$ coprime with $l$, $\,{\rm H}_i (Z_0;\,\Z / k\Z ) =
0,\,i=1,\dots,n$, i.e. the fibre $Z_0$ is $\Z / k\Z$--acyclic.

\noin Then the zero fibre $X_0 = p^{-1} (0) \subset \C^{n+1}$ of the polynomial
$p(x,\,z) = q(x) + z^k$ is diffeomorphic to $\R^{2n}$.} \\

\noin Note that the polynomial $p$ is a quasi--invariant of weight $kl$ of the
$\C^*$--action ${\tilde G}(\lambda,\,x,\,z) = (G(\lambda^k,\,x),\,\lambda^l z)$
on $\C^{n+1}$. The morphism $\pi\,:\,X_0 \ni(x,\,z) \longmapsto x \in \C^n$
represents $X_0$ as a $k$--fold branched cyclic covering over $\C^n$  ramified
along $Z_0$. This covering  is equivariant with respect to the actions $G$ on
$\C^n$ and ${\tilde G}\,|\,X_0$ on $X_0$.

The generic fibre $X_1 = p^{-1} (1)$ topologically is the joint $Z_1 \ast \Z /
k\Z$, where $Z_1 = q^{-1} (0) \subset \C^n$ is the generic fibre of $q$ [Ne].
Therefore, $X_1$ is not contractible as soon as $Z_1$ is not contractible [Ka
1, Lemma 9].

The assumption (i) is always fulfilled for a generic fibre of a primitive
polynomial [Ka 1, Lemma 8]. Though the zero fibre of a $\C^*$--quasi--invariant
is usually non--generic, there exist non--trivial examples where this and all
the other conditions of Theorem 4.3 are satisfied. \footnote{Sh. Kaliman has
informed me that in fact the condition (i) of Theorem 4.3 is superfluous; see
Appendix.} \\

\noin {\bf 4.4. Proposition} [Ka 1, Theorem 10]. {\it Put $$q(x,\,y,\,z) = x +
x^ay^b + y^cz^d\,,$$  $l = bd$ and $G(\lambda,\,x,\,y,\,z) = (\lambda^l
x,\,\lambda^{-r} y, \,\lambda^s z)$, where $r = d(a-1),\, s = c(a-1) + b$. If
$(s,\, d) = 1$, then the polinomial $q$ and the $\C^*$--action $G$ verify all
the assumptions of Theorem 4.3. Furthermore, the Euler characteristic of the
generic fibre $Z_1 = q^{-1} (1)$ is different from $1$ as soon as $b,\,d \ge
2$. } \\

\noin Let $p(x,\,y,\,z,\,t) = q(x,\,y,\,z) + t^k$, where $q$ is as above and
$k$ is a prime such that $(bd,\,k) = 1$. Then, by Theorem 4.3, $X_0 =
X_{a,\,b,\,c,\,d,\,k} = p^{-1} (0)$ is a smooth contractible hypersurface in
$\C^4$. Later on we will see that most of these threefolds are exotic
$\C^3$--s.

The polynomial $p$ being $\tilde G$--quasi--invariant, the threefold $X_0$
carries a $\C^*$--action ${\tilde G}\,|\,X_0$. In general, a polynomial $f$ on
$\C^n$ is a quasi--invariant of a linear diagonalized $\C^*$--action iff its
Newton diagram is linearly degenerate, i.e. if it lies in an affine hyperplane.
This is always the case when $f$ consists of $n$ monomials only, like in the
preceeding examples. Since all the hypersurfaces discussed in this section
carry $\C^*$--actions, none of them is of log--general type. However, as we
will see, for some of them the log--Kodaira dimension is non--negative.

The next result provides an estimate from below of the log-Kodaira dimensions
of ramified coverings, and so makes it possible, in certain cases, to recognize
exotic $\C^n$--s among hypersurfaces as in Theorem 4.3. \\

\noin {\bf 4.5. Proposition} [Ka 1, Corollaries 12, 13]. {\it Let $f\,:\,W \to
V$ be a morphism  of smooth quasiprojective varieties, which is a branched
covering ramified over a divisor $R \subset V$ of simple normal crossing type.
Assume that the Sakai analytic dimension $k_c$ of the complement $V \setminus
R$ is non--negative  (see [Sa]) \footnote {being non--negative $k_c$ must
coincide with the log--Kodaira dimension $\k$ [Sa]}. If the ramification order
of $f$ on each of the irreducible components of $f^{-1} (R) \subset W$ is high
enough, then $\k (W) \ge \k (V \setminus R)$. Consequently, in this case $W$ is
of log-general type if $V \setminus R$ is so.} \\

More carefull analysis with the same type of arguments [Ka 1, Theorem B] leads
to the conclusion that $\k (X_{a,\,b,\,c,\,d,\,k}) = 2$ if $k > a \ge 4, \,d=
a-1, \,(b, \,d) = (bd, \,k) = 1$ and $ (b,\,c) > d^2 k$. For instance, $\k(X) =
2$ for $X = X_{4,\,46,\,92,\,3,\,5} = \{x + x^4 y^{46} + y^{92}z^3 + t^5 = 0\}
\subset \C^4$. See [Ru 1, Ka 3] for some other cases when the log--Kodaira
dimension of these threefolds is at least non--negative or more. Note that an
exotic threefold of non--negative log--Kodaira dimension cannot be an exotic
product structure on $\C^3$, which is always of log--Kodaira dimension
$-\infty$. \\

\centerline {\bf 5. Russell's $\C^*$--threefolds and the Makar--Limanov
invariant} \\

By {\it a $\C^*$--variety} we mean a smooth irreducible algebraic variety
endowed with a regular effective $\C^*$--action. Most of the exotic $\C^n$--s
which are known, except some Kaliman modifications, are $\C^*$--varieties.
Thus, we come to the following \\

\noin {\bf 5.1. Problem.} {\it Classify contractible $\C^*$--varieties up to
equivariant isomorphism.} \\

This includs the famous linearization problem for $\C^*$--actions on the affine
space $\C^3$, which is still open (see e.g. [KoRu, Kr 2]). \\

\noin {\bf 5.2. Koras--Russell bicyclic covering construction.} Analysing
Dimca's examples 4.2 from the point of view of the previous work with M. Koras
[KoRu], P. Russell [Ru 1] came to a remarkable general method of constructing
contractible $\C^*$--threefolds. In particular, it yields all of them of a
certain restricted type (namely, {\it tame of mixed type}; see (5.4) and
Theorem 5.5 below), including those of 4.2 and 4.4 above. We discuss here some
principal points of this construction.

Denote by $\omega_r$ the cyclic group of the complex $r$--roots of unity. Let
$B$ be a smooth contractible algebraic variety and $Z_1,\, Z_2 \subset B$ be
two smooth divisors which meet normally. For a pair $(\alpha_1,\,\alpha_2)$ of
coprime positive integers consider the bicyclic covering over $B$ branched to
order $ \alpha_i$ over $Z_i, \,i=1,2$. We get a commutative diagram \\

\begin{picture}(200,100)
\unitlength0.2em
\thicklines
\put(75,39){$X$}
\put(51,18){$X_1$}
\put(97,18){$X_2$}
\put(75,-2){$B$}

\put(72,37){$\vector(-1,-1){15}$}
\put(80,37){$\vector(1,-1){15}$}
\put(58,16){$\vector(1,-1){15}$}
\put(96,16){$\vector(-1,-1){15}$}
\put(77,34){$\vector(0,-1){28}$}

\put(57,32){$\omega_{\alpha_2}$}
\put(89,32){$\omega_{\alpha_1}$}
\put(57,6){$\omega_{\alpha_1}$}
\put(90,6){$\omega_{\alpha_2}$}
\put(79,20){$\omega_{\alpha_1 \alpha_2}$}

\put(160,22){$(2)$}
\end{picture}\\

\noin {\bf Question.} {\it Under which assumptions on $Z_1,\,Z_2$ the resulting
variety $X$ is contractible?} \\

\noin Here $X$ is acyclic if both $\omega_{\alpha_1}, \, \omega_{\alpha_2}$ act
trivially on the homologies $H_* (X; \Z)$. Indeed, in this case the homologies
would be $\alpha_i$--torsions, $i = 1,2$, and hence trivial [Ru 1]. This is so
if $B$ is a $\C^*$--variety and the divisors $Z_1,\, Z_2$ are invariant under
the $\C^*$--action and, furthermore, they are given as $Z_i = q_i^{-1}
(0),\,i=1,2$, where $q_i \in \C[B]$ is a regular $\C^*$--quasi--invariant of
weight which is relatively prime with $\alpha_i, \,i=1,2$. Indeed, under these
assumptions (2) is a diagram of equivariant morphisms of $\C^*$--varieties, and
$\omega_{\alpha_i}, i=1,2$, acts on $X$ via the $\C^*$--action. Therefore, it
induces the trivial action in the homologies. We also need, of course, $X$ to
be simply connected (cf. Appendix).

{}From now on we restrict the consideration by $\C^*$--threefolds, and namely
by those of mixed (or hyperbolic) type. We say that a $\C^*$--threefold $X$
with a $\C^*$--action $G$ is {\it of mixed type} if $G$ has the unique fixed
point $X^G = \{x_0\}$ and the weights of the diagonalized linear action
$dG\,|\,T_{x_0} X$ are of different signs. One may assume that $$dG
(\lambda,\,x,\,y,\,z) = (\lambda^{-a} x,\, \lambda^b y,\,\lambda^c
z),\,\,\lambda \in \C^* ,$$
\noin where $a,\,b,\,c > 0$. Since $G$ is effective, $a,\,b,\,c$ are relatively
prime: $(a,\,b,\,c) = 1$. The triple $(a,\,b,\,c)$ is called {\it reduced} if,
moreover, $(a,\,b) = (b,\,c) = (a,\,c) = 1$.

We will say that a mixed $\C^*$--threefold $X$ is {\it tame} if the algebraic
quotient $X//G$ is isomorphic to $\C^3 //dG$, or, what is the same, to one of
the surfaces $\C^2 // \omega_r$, where $\omega_r$ acts diagonally in $\C^2$.
$X$ is said to be {\it linearizable} if it is equivariantly isomorphic to
$\C^3$ with a linear $\C^*$--action. The following statement is unexplicitly
contained in [KoRu, Ko]. \\

\noin {\bf 5.3. Proposition.} {\it If $X$ is a tame contractible affine
$\C^*$--threefold of mixed type with a reduced triple of weights at the fixed
point, then $X$ is linearizable.} \\

The next example due to P. Russell shows the importance of the assumption on
$X$ being tame. \\

\noin {\bf Example.} Let $S$ be a contractible surface of non--negative
log--Kodaira dimension. Put $Y = S \times \C$, and let $X$ be the Kaliman
modification of $Y$ along $S \times \{0\}$ with center at a point $\{s_0\}
\times \{0\}$. The tautological $\C^*$--action on $\C$ lifts to a
$\C^*$--action on $Y$, which in turn induces a $\C^*$--action of mixed type on
$X$ with the reduced triple of weights $(-1,\,-1,\,1)$. Here $X// \C^* \simeq S
\not\simeq \C^2 // \omega_r$, and being an exotic $\C^3$, $\,X$ is not
linearizable. \\

Now the idea of [KoRu, Ru 1] can be described as follows. Starting with a tame
contractible $\C^*$--threefold $X$ with a non--reduced triple of weights
$(-a,\,b,\,c)$, we may factorise it by a cyclic subgroup of $\C^*$ to make the
triple of weights being reduced. By Proposition 5.3, we arrive in this way to
$\C^3$ with a linear $\C^*$--action. Namely, put $$\alpha = (b,\,c),\,\beta =
(a,\,c),\,\gamma = (a,\,b)\,\,\,\,{\rm and}\,\,\,\,a' = a/\beta\gamma,\,b' =
    b/\alpha\gamma,\,c' = c/\alpha\beta\,.$$

\noin Then both $(\alpha,\,\beta,\,\gamma)$ and $(a',\,b'\,c')$ are reduced
triples and $(a',\,\alpha) = (b',\,\beta) = (c',\,\gamma) = 1$. Furthermore, $B
= X / \omega_{\alpha\beta\gamma}$ is a tame $\C^*$--threefold of mixed type
with the reduced triple of weights $(-a',\,b',\,c')$, which is, due to
Proposition 5.3, isomorphic to $\C^3$ with a linear $\C^*$--action of the same
type. \\

\noin {\bf 5.4. Russell's threefolds.} P. Russell [Ru 1] realized the converse
procedure. Starting now with $\C^3$ with a linear diagonalized reduced mixed
$\C^*$--action and passing to the corresponding tricyclic coverings, he
reconstructed all possible tame contractible affine $\C^*$--threefolds of mixed
type. In what follows we call them {\it Russell's threefolds}.

The principal point of Russell's construction is the choice of branching
divisors $Z_0,\,Z_1,\,Z_2 \subset \C^3$ of the tricyclic covering. The first of
them $Z_0$ appears naturally. Indeed, let $X$ be a contractible
$\C^*$--threefold of the mixed type $(-a,\,b,\,c)$ with the fixed point $x_0
\in X$. Put $$X^+ = \{x \in X\,|\,\lim_{\lambda \to 0} G_{\lambda} x = x_0
\}\,,$$ $$X^- = \{x \in X\,|\,\lim_{\lambda \to 0} G_{\lambda^{-1}} x = x_0
\}\,.$$

\noin Then $X^+ ,\,X^-$ are isomorphic to $\C^2$ and $\C$ respectively and meet
transversally. We call them {\it the positive plane} resp. {\it the negative
axis}. If $\sigma\,:\,X \to B,\,B = X / \omega_{\alpha\beta\gamma} \simeq
\C^3$, is the quotient morphism, then $\sigma (X^{\pm}) = B^{\pm}$ are,
respectively, the coordinate plane $x = 0$ and the coordinate axis $y = z = 0$,
and $\sigma$ is branched to order $\alpha$ along $Z_0 = B^+$.

The two other divisors $Z_1,\,Z_2$ are chosen as follows. In order to get $X_i
\simeq \C^3$ in diagram (2) we take $Z_i \simeq \C^2,\,i=1,2$. Moreover, after
passing to the first covering we would like to have in $X_i \simeq
\C^3\,(i=1,2)$ the situation described in Theorem 4.3. To this point these
three embedded planes $Z_0,\,Z_1,\,Z_2 \subset \C^3$ should satisfy the
following conditions:\\{\it

\noin i)  $Z_i$ is equivalent  to a coordinate plane under a tame automorphism
of $\C^3,\, i=1,2$;

\noin ii) $Z_0 \cup Z_1 \cup Z_2$ is a normal crossing divisor;

\noin iii) $Z_i$ are invariant under the $\C^*$--action $(\lambda,\,x,\,y,\,z)
\longmapsto (\lambda^{-a'} x,\,\lambda^{b'} y,\,\lambda^{c'} z)$ on $\,\C^3,\,
i=0,1,2$;

\noin iv) locally at the origin the triple $Z_0,\,Z_1,\,Z_2$ determines a
quasi--homogeneous coordinate system \footnote{this follows, of course, from
ii) and iii).}

\noin v) the intersection $Z_1 \cap Z_2$ consists of the negative axis $B^-$
and of $r-1$ closed $\C^*$--orbits, where $r \ge 1$.} \\

\noin If $Z_1,\,Z_2$ can be linearized simulteneously (this corresponds to $r =
1$), then clearly $X \simeq \C^3$. Therefore, in interesting cases $r > 1$.

Put $C_i = Z_i \cap P,\,i=1,2$, where $P = \{x = 1\} \subset \C^3$. Note that
the surface $Z_i \subset \C^3$ is the closure of the orbit of the curve $C_i$
under the $\C^*$--action. The affine plane $P$ is invariant with respect to the
induced $\omega_{a'}$--action, and by (iii) the curves $C_1, \,C_2 \subset P$
should be $\omega_{a'}$--invariant, too. Furthermore, $C_1,\,C_2$ are
isomorphic to $\C$ and meet normally at the origin $P \cap B^-$ and in $r-1$
other points. If $C_i$ are given in the plane $P$ by the equations $p_i(y,\,z)
= 0, \,p_i \in \C[y,\,z],\,i=1,2$, then the equations of $Z_i$ are $\p_i
(x,\,y,\,z) = 0$, where $\p_i \in \C [x,\,y,\,z], \,i=1,2$, are defined as
follows: $$\p_1 (x^{a'},\,y,\,z) = x^{-b'} p_1 (x^{b'} y,\,x^{c'} z)\,\,\,{\rm
and}\,\,\,\p_2 (x^{a'},\,y,\,z) = x^{-c'} p_2 (x^{b'} y,\,x^{c'}
z)\,\,\,\,\,\,\,\,\,\,\,\,\,\,\,(3)$$ As soon as a pair of plane curves
$C_1,\,C_2$ as above is chosen in such a way that $\p_i$ are polynomials, $i =
1,2$, the corresponding triple $Z_0,\,Z_
1,\,Z_2$ satisfies all the conditions (i)-(v). \\

\noin {\bf 5.5. Theorem.} a) [KoRu, Ru 1] {\it Fix two reduced triples
$(a',\,b',\,c')$ and $(\alpha,\,\beta,\,\gamma)$ of positive integers such that
$(a',\,\alpha) = (b',\,\beta) = (c',\,\gamma) = 1$. Let $C_i = p_i^{-1} (0),
\,Z_i = \p_i^{-1} (0),\,i=1,2,$ and $Z_0 = \{x = 0\}$ be as above. Let $X \to
\C^3$ be the tricyclic covering ramified to order $\alpha$ over $Z_0$, to order
$\beta$ over $Z_1$ and to order $\gamma$ over $Z_2$\footnote{that is, $X = {\rm
spec}\,A$, where $A$ is the extension of $\C[x,\,y,\,z]$ by the corresponding
roots of $x,\,\p_1$ and $\p_2$.}. Then $X$ is a Russell threefold. Conversely,
any Russell threefold is obtained by the above construction.}

\noin b) [ML, KaML 1,2] {\it X as above is an exotic $\C^3$ except in the cases
when $$(r-1)(\beta - 1)(\gamma - 1) = 0\,.$$} {\bf 5.6. Remark.} Putting $p_2
(y,\,z) = z$, which is possible, we may present a Russell threefold $X$ as the
hypersurface in $\C^4$ with the equation $$\p_1(x^{\alpha},\,y,\,z^{\gamma}) +
t^{\beta} = 0\,.$$

\noin {\bf 5.7. Examples.} Put $(-a', \,b',\,c') = (-r + 1,\,1,\,1)$, where
$r \ge 2$, and $p_1 (y,\,z) = y + y^r + z$. Then $\p_1 (x,\,y,\,z) = y + x y^r
+ z$ and $$X = \{y + x^{\alpha} y^r + z^{\gamma} + t^{\beta} = 0\} \subset
\C^4\,.$$ In the simplest non--trivial case $r=2,\,\alpha = 1,\,\beta = 3,
\gamma = 2$ we get the affine cubic $X_0 \subset \C^4$ with the equation $$y +
xy^2 + z^2 + t^3 = 0\,,$$ which is an exotic $\C^3$ [ML]. \\

The following result makes precise the statement of Theorem 5.5, b). It was
obtained by a rather elementary method of analyzing the defining equations of
Russell's threefolds [KaML 1,3],  \\

\noin {\bf 5.8. Theorem} [KaML 1]. {\it Let $X$ be a Russell threefold
constructed by the data $(r,\,\alpha,\,\beta,\,\gamma)$, where $r > 1,\,\alpha
\ge 2,\, \beta,\,\gamma \ge 4$. Then there is no dominant morphism $\C^3 \to
X$. }\\

It was previously known [Ru 1] that $\k (X) = 2$ if $ \alpha \ge a'\beta\gamma$
and $\beta,\,\gamma >> 1$; furthermore, $\k (X) \ge 0$ if  $a' = 1,\, \alpha
\ge 2, \, \beta,\,\gamma \ge 4$. But $\k (X) = -\infty$ for the Russell's
threefolds $X$ as in 5.7 with $\alpha = 1$, since the complement $X \setminus
\{y = 0\}$ is isomorphic to $\C^* \times \C^2$. Furthermore, the cubic $X_0$
from 5.7 is dominated by $\C^3$. \\

\noin {\bf 5.9. The Makar--Limanov invariant} \\

\noin {\bf Definition.} Recall that a derivation $\partial$ of a ring $A$ is
called {\it locally nilpotent} if each element $a \in A$ is vanished by an
appropriate positive power $n = n(a)$ of $\partial$, i.e. $\partial^n (a) = 0$.
Denote $A^{\partial} = {\rm Ker}\,\partial$; $A^{\partial}$ is called {\it the
ring of constants of} $\partial$. Let $LN(A)$ denote the set of all locally
nilpotent derivations of $A$. Put $A_0 = \bigcap_{\partial \in LN(A)}
A^{\partial}$. We call $A_0$ {\it the ring of absolute constants of $A$}. Note
that $A_0 = \{\C\}$ if $A = \C[x_1,\,\dots,\,x_n]$ is a polynomial ring. \\

The subring $A_0 \subset A$ of the absolute constants is invariant under ring
isomorphisms. It was introduced in [ML], where it was shown that it is
non--trivial in the case of the algebra $A=\C[X_0]$ of the regular functions on
the Russell cubic $X_0$ (see 5.7). We call $A_0$ {\it the Makar--Limanov
invariant of $A$}. \\

\noin {\bf 5.10. Theorem} [KaML 2, Theorem 8.3]. {\it Let $A = \C[X]$, where
$X$ is a Russell threefold. Then $A_0 = A$ (i.e. there is no locally nilpotent
derivations on $A$) except in the following two cases:

\noin i) if $X = \{x + x^ry + z^{\beta} + t^{\gamma} = 0\}$, then $A_0 =
\C[x]$;

\noin ii) if $X = \{x + (x^r + z^{\beta})^l y + t^{\gamma} = 0\}$, then $A_0 =
\C[x,\,z]$.} \\

This proves (b) of Theorem 5.5. \\

\vs

\centerline {\bf 6. APPENDIX: Simply connectedness of $\C^*$--equivariant
cyclic coverings}\\

The results of this section are due to Sh. Kaliman \footnote{Letter to the
author from 12.03.1995. We get them placed here with the kind permission of Sh.
Kaliman and P. Russell.}. The original presentation has been modified by P.
Russell by picking out the group theoretic component (see Proposition 6.2). In
particular, he used the following  \\

\noin {\bf 6.1. Definition}. Let $G$ be a group. We say that a subgroup $H
\subset G$ is {\it normally generated by elements} $a_1,\dots, a_n \in H$ if it
is generated by the set of all elements conjugate to $a_1,\dots, a_n$. Thus,
$H$ is the minimal normal subgroup of $G$ that contains $a_1,\dots, a_n$. We
denote it by $\,<<a_1,\dots, a_n>>$.

We say that $G$ is  {\it normally one--generated} if $G = \,<<a>>\,$ for some
element $a \in G$. If $A, \,B \subset G$, then $[A, \,B]$ denotes the subgroup
generated by all the commutators $[a,\,b] = aba^{-1}b^{-1}$, where $a \in A,\,b
\in B$. \\

\noin {\bf 6.2. Proposition.} {\it Put $K = [G,\,G]$. Assume that $G =
\,<<a>>\,$ is normally one--generated and that $G_{\rm ab} \simeq \Z$ where
$G_{\rm ab} = G/K$ is the abelianization of $G$. Then the following statements
hold.\\

\noin a) $K = \,<<[a,\,K]>>$.\\

\noin b) Suppose that $[a^l,\,K] = {\bf 1}$ for some $l \neq 0$. Fix any $k \in
\Z$ with $(k,\,l) = 1$. Let $G_k = \rho^{-1} (H_k)$, where $\rho\,:\,G \to
G_{\rm ab} = G/K$ is the canonical epimorphism and $H_k \subset G_{\rm ab}$ is
the subgroup generated by $a^k ({\rm mod}\,K)$. Then $G_k = \,<<a^k>>$. }\\

\noin {\it Proof.} a) Denote by $A$ the cyclic subgroup generated by $a$: $A =
\,<a>\, \subset G$. Set $M_a = \,<< \,[A,\,G]\, >>\, = \,<<\,[a^k,\,G]\,|\,k
\in \Z\,>>$. \\

\noin {\bf Claim 1.} $K = M_a$. \\

\noin {\it Proof.} The abelianization $G_{\rm ab}$ of $G$ is a free cyclic
group generated by the class $a \,({\rm mod}\,K) \neq 0$. From the short exact
sequence $${\bf 1} \to K \to G \to G_{\rm ab} \simeq \Z \to {\bf 0}$$ it
follows that the commutator subgroup $K = [G,\,G]$ consists of the elements $$g
= \prod\limits_{i=1}^r c_ia^{m_i} c_i^{-1} \in G$$ such that
$\sum\limits_{i=1}^r m_i = 0$.

Set $k_0 = 0$ and $k_i = \sum\limits_{j=1}^i m_j$, so that $m_i = -k_{i-1} +
k_i$. Denote by $\sim b$ any element conjugate to $b$. Since $ca^{k+l}c^{-1} =
(ca^kc^{-1})(ca^lc^{-1})$, with the above notation every element $g \in G$ can
be written as \footnote{Although the equality $\sim a^{k + l} = \,(\sim a^k) (
\sim a^l)\,\,\,$ is not symmetric any more, this does not cause problems, as
well as the use of the non--commutative product symbol $\prod$.} $$g =
\prod\limits_{i=1}^r (\sim a^{m_i} ) = \prod\limits_{i=1}^r (\sim a^{-k_{i-1} +
k_i}) = (\prod\limits_{i=1}^{r-1} (\sim a^{k_i})(\sim a^{-k_i}))\,a^{k_r}\,,$$
where $k_r = 0$ iff $g \in K$. Furthermore, note that $$(\sim a^k) ( \sim
a^{-k}) = ca^kc^{-1}da^{-k}d^{-1} = c(a^kc^{-1}da^{-k}d^{-1}c)c^{-1} =
c[a^k,\,b]c^{-1} = \,\sim [a^k,\,b]\,,$$ where $b = c^{-1}d$. Thus, we have $$g
= (\prod\limits_{i=1}^{r-1} (\sim[a^{k_i} ,\,b_i]))\,a^{k_r}\,.$$ If $g \in K$,
then $k_r = 0$ and, therefore, $g \in M_a$. This proves the inclusion $K
\subset M_a$. Since, evi

dently, $M_a \subset K$, we have $K = M_a$. \qed \\

Put $$N_a = \,<<[a^{\epsilon},\,G]\,|\,\epsilon = \pm1>>\, \,\subset M_a\,.$$

\noin {\bf Claim 2.}  $M_a = N_a$. \\

\noin {\it Proof.} From the identity $$[a^k,\,b] =
(a^{k-1}[a,\,b]a^{-k+1})[a^{k-1},\,b]\,,\,\,\,k\in\Z_{>0}\,,$$ we obtain by
induction $$[a^k,\,b] =
(a^{k-1}[a,\,b]a^{-k+1})(a^{k-2}[a,\,b]a^{-k+2})\dots(a[a,\,b]a^{-1})[a,\,b]
\in N_a\,.$$ To show that $[a^{-k} ,\,b] \in N_a,\,k\in\Z_{>0}$, it is enough
to replace $a$ by $a^{-1}$ in the above identities. Thus, $M_a \subset N_a$.
\qed \\

\noin {\bf Claim 3.} $N_a = <<[a,\,K]>>$. \\

\noin {\it Proof.} Since $G_{\rm ab} = <a ({\rm mod}\,K) >$ is a cyclic group,
any element $b \in G$ can be written as $b = a^md$, where $d \in K$. Therefore,
$$[a^{\epsilon},\,b] = a^{\epsilon}(a^m d)a^{-\epsilon}(a^m d)^{-1} = a^m
(a^{\epsilon} d a^{-\epsilon} d^{-1} )a^{-m} =\, \sim [a^{\epsilon} ,\,d] \in
\,<<[a^{\epsilon} ,\,K]>>\,.$$ Finally, $$[a^{-1},\,b]^{-1} = a^{-1}[a,\,b] a
\in \,<<[a,\,K]>>\,,$$ and hence $[a^{-1},\,b] \in \,<<[a,\,K]>>$. \qed\\

Now (a) follows from Claims 1--3. \qed\\

To prove (b), we start with the following \\

\noin {\bf Claim 4.} {\it Under assumptions of (b), $K = \,<<[a^k,\,K]>>$.} \\

\noin {\it Proof.} Represent $1 = \mu l + \nu k$, where $\mu, \nu \in \Z$. Then
for $d \in K$ we have $$[a^{\epsilon} ,\,d] = a^{\epsilon} d a^{-\epsilon}
d^{-1} = (a^k)^{\epsilon\nu} d (a^k)^{-\epsilon\nu} d^{-1} \in
\,<<[a^k,\,K]>>\,.$$ Since $\,<<[a^k,\,K]>>\, \subset K$ and by (a) $K =
\,<<[a,\,K]>>\,$, the Claim follows. \qed \\

Furthermore, note that for any $c \in G_k$ we have $c = a^{mk} d$, where $d \in
K$. Therefore, to prove (b) it suffices to show that $K \subset \,<<a^k>>$.

Take $[g,\,h] \in K$ arbitrary. Due to Claim 4 we have the presentation
$$[g,\,h] = \prod\limits_{i=1}^N d_i[a^{\nu_i k},\,c_i] d_i^{-1} =
\prod\limits_{i=1}^N (d_i a^{\nu_i k}d_i^{-1})((d_ic_i)a^{-\nu_i
k}(d_ic_i)^{-1}) $$ $$= \prod\limits_{i=1}^N (\sim a^{\nu_i k}) (\sim a^{-\nu_i
k}) \in \,\, <<a^k>>\,.$$ This proves (b). \qed \\

\noin {\bf 6.3. Notation.} Let $X$ be a smooth irreducible algebraic variety,
$q \in \C[X],\,F_0 = q^{*}(0)$ and $F_1 = q^{-1}(1)$. Assume that $F_0$ is a
reduced and irreducible divisor. Fix a smooth complex disc $\Delta\subset X$
which meets $F_0$ normally at a smooth point of $F_0$, and a small positive
simple loop $\delta \subset \Delta$ around $F_0$. It defines uniquely up to
conjugacy an element $\alpha \in \pi_1(X \setminus F_0)$ \footnote{indeed, this
easily follows from the connectedness of the smooth part ${\rm reg}\,F_0$ of
$F_0$.}. Following Fujita [Fu 2, (4.17)] we call such an $\alpha$ {\it the
vanishing loop of the divisor $F_0$}, and the group $\,<<\alpha>>\, \subset
\pi_1(X \setminus F_0)$  {\it the vanishing subgroup of $F_0$}, keeping in mind
that $\,<<\alpha>>\,$ is contained in the kernel of the natural surjection
$i_*\,:\,\pi_1(X \setminus F_0) \to \pi_1(X)$.\\

The following statement should be well known. However, in view of the lack of
references we sketch the proof. \\

\noin {\bf 6.4. Lemma.} ${\rm Ker}\,i_* = \,<<\alpha>>\,$. \\

\noin {\it Proof.} Let a loop $\gamma\,:\,S^1 \to X \setminus F_0$ represents
the class $[\gamma] \in {\rm Ker}\,i_*$, i.e. $\gamma$ is contractible in $X$.
Fix a stratification of $F_0$ which satisfies the Whitney condition A and
contains the regular part ${\rm reg}\,F_0$ of $F_0$ as an open stratum. By the
Thom Transversality Theorem the homotopy $S^1 \times [0,\,1] \to X$ of $\gamma
= \gamma_0$ to the constant loop $\gamma_1 \equiv {\rm const}$ can be chosen
being transversal to the stratification, and therefore such that its image
meets the divisor $F_0$ in a finite number of its regular points only. We may
also assume that these intersection points $p_1,\,\dots,\,p_n \in {\rm
reg}\,F_0\,, \,\,\,p_i \in \gamma_{t_i} \cap F_0$, correspond to different
values $0 < t_1 < \dots < t_n < 1$ of the parameter of homotopy $t \in [0, 1]$.
If $s_i \in [0,\,1],\,\,s_i < t_i < s_{i+1},$ and ${\bar \gamma}_i =
\gamma_{s_i}\,:\,S^1 \to X \setminus F_0 , \, i=1,\dots,n+1$, then clearly
${\bar \gamma}_{i+1}^{-1}\cdot {\b
ar \gamma}_i \approx \delta_i^{\epsilon_i}$, i.e. ${\bar \gamma}_i \approx
{\bar \gamma}_{i+1}\cdot\delta_i^{\epsilon_i}$, where $\delta_i$ is a vanishing
loop of $F_0$ at the point $p_i$ and $\epsilon_i = \pm 1$, and ${\bar
\gamma}_{n+1} \approx {\rm const}$. Thus, $[\gamma] = [{\bar \gamma}_1] =
[\delta_n]^{\epsilon_n}\cdot \dots \cdot [\delta_1]^{\epsilon_1} \in
\,<<\alpha>>\,$, and we are done. \qed\\

\noin {\bf 6.5. Corollary.} {\it If $\pi_1 (X) = {\bf 1}$, then the group $G =
\pi_1 (X \setminus F_0)$ is normally one--generated by the vanishing loop
$\alpha$ of $F_0$.} \\

\noin {\bf 6.6.} Assume further that the restriction $q\,|\,(X \setminus
F_0)\,:\,X \setminus F_0 \to \C^*$ is a smooth fibration. Then we have the
exact sequence
$$\,\,\,\,\,\,\,\,\,\,\,\,\,\,\,\,\,\,\,\,\,\,\,\,\,\,\,\,\,\,\,\,\,\,\,\,\,\,\,\,\,\,\,\,{\bf 1} \to \pi_1(F_1) {\stackrel{i_*}{\longrightarrow}} \pi_1 (X \setminus F_0) {\stackrel{q_*}{\longrightarrow}} \Z \to {\bf 0} \,\,\,\,\,\,\,\,\,\,\,\,\,\,\,\,\,\,\,\,\,\, \,\,\,\,\,\,\,\,\,\,\,\,\,\,\,\,\,\,\,\,\,\,\,\,\,\,\,\,\,\,\,\,\,\,\,\,\,\,\,\,\,\,\,\,(4)$$ such that $q_*(\alpha) = 1 \in \Z$.  \\

\noin {\bf 6.7. Lemma.} {\it In the assumptions as above suppose additionally
that $X$ is simply connected. Then $i_*\pi_1(F_1) = K = [G,\,G]$, and $G_{\rm
ab} = H^1 (X \setminus F_0) \simeq \Z$.} \\

\noin {\it Proof.} By Corollary 6.5 we have $G = \,<<\alpha>>\,$. From the
exact sequence (4) it follows that $$i_*\pi_1(F_1) = {\rm Ker}\,q_* = \{g =
\prod\limits_{i=1}^r (\sim\alpha^{k_i})\,\,|\,\,\sum\limits_{i=1}^r k_i = 0\} =
K$$ (see the proof of Claim 1 in Proposition 6.2). This proves the first
assertion. The second one follows by applying (4) once again. \qed\\

\noin {\bf 6.8. Lemma.} {\it Let, in the notation as above, $q$ be a
quasi--invariant of a positive weight $l$ of a regular $\C^*$--action $t$ on $X
\setminus F_0$. Then $[\alpha^l,\,K] = {\bf 1}$.} \\

\noin {\it Proof.} Let $\varphi_l\,:\,Y_l \to X \setminus F_0$ be the cyclic
covering of order $l$: $$Y_l = \{(x,\,z) \in (X \setminus F_0) \times
\C\,|\,z^l = q(x)\}\,.$$ Put $q_l = q \circ \varphi_l\,:\,Y_l \to \C^*$. Define
a morphism $\theta\,:\,F_1 \times \C^* \to Y_l$ as follows: $$\theta
(x,\,\lambda) = (t(\lambda,\,x),\,\lambda)\,, \,x\in F_1,\,\lambda \in \C^*\,.
$$ It is easily seen that $\theta\,:\,F_1 \times \C^* \to Y_l$ is an
isomorphism. We have $$\varphi_* \pi_1(Y_l) = \varphi_* (\pi_1(F_1) \times \Z)
= <i_*\pi_1(F_1),\,\alpha^l>\,.$$ This implies that $\alpha^l$ commutes with $K
= i_*\pi_1(F_1)$. \qed \\

The next theorem is the main result of this section. \\

\noin {\bf 6.9. Theorem (Sh. Kaliman).} {\it Let $X$ be a simply connected
smooth irreducible algebraic variety, $q \in \C[X]$ be a regular function on
$X$ such that \\

\noin i) $F_0 = q^*(0)$ is a smooth reduced irreducible divisor, and \\

\noin ii) $q\,|\,(X \setminus F_0)$ is a quasi--invariant of weight $l > 0$ of
a regular $\C^*$--action $t$ on $X \setminus F_0$. \\

\noin Let $\sigma_k \,:\,X_k \to X$ be the cyclic covering branched to order
$k$ over $F_0$: $$X_k = \{(x,\,z) \in X \times \C\,|\,q(x) =
z^k\}\,,\,\,\,\,\,\,\,\,\sigma_k (x,\,z) = x\,.$$ If $(k,\,l) = 1$, then $X_k$
is simply connected. } \\

\noin {\it Proof.}  Put $q_k = q\circ \sigma_k \in \C[X_k]$ and $F_{k,0} =
q_k^{-1}(0) \subset X_k$. Since $X_k \setminus F_{k,0} \to X \setminus F_0$ is
a non--ramified $k$-sheeted cyclic covering, the induced homomorphism
$$(\sigma_k)_*\,:\,\pi_1(X_k \setminus F_{k,0}) \to \pi_1(X \setminus F_0) =
G$$ is an injection onto a normal subgroup $G_k$ of $G$ of index $k$, and
$G/G_k \simeq \Z/k\Z$. Clearly, $\alpha^k \in G_k$ is covered by a vanishing
loop $\beta \in \pi_1(X_k \setminus F_{k,0})$ of the smooth divisor $F_{k,0}
\subset X_k$, i.e. $(\sigma_k)_* (\beta) = \alpha^k$. Therefore,
$\,<<\alpha^k>> \,\subset G_k$.

In fact, $G_k$ has the same meaning that in Proposition 6.2(b), i.e. $G_k =
q_*^{-1} (H_k)$, where $H_k = k\Z \subset \Z \simeq G_{\rm ab}$. Indeed, by the
universal property of the commutator subgroup, under the homomorphism
$\tau\,:\, G \to G/G_k \simeq \Z/k\Z$ we have $K \subset {\rm Ker}\,\tau =
G_k$, and hence $G_k = q_*^{-1} (q_* (G_k))$. Furthermore, $q_* (G_k) \supset
k\Z = H_k$, because $\alpha^k \in G_k$ and $q_* (\alpha_k)=1 \in \Z$. Actually,
 $q_* (G_k) = H_k$, since $[G\,:\,G_k] = k$. It follows that $G_k = q_*^{-1}
(H_k)$.

By Lemma 6.8, we have $[\alpha^l,\,K] = {\bf 1}$, so that Proposition 6.2(b)
can be applied. Due to this Proposition, $G_k = \,<<\alpha^k>>\,$. Or, what is
the same, $\pi_1(X_k \setminus F_{k,0}) = \,<<\beta>>\,$. The inclusion
$i\,:\,X_k \setminus F_{k,0} \hookrightarrow X_k$ induces an epimorphism
$i_*\,:\,\pi_1(X_k \setminus F_{k,0}) \to \pi_1(X_k)$ with the kernel
$\,<<\beta>>\,$ (see Lemma 6.4). Thus, $\pi_1(X_k) = {\bf 1}$, as desired. \qed

\vs

\centerline {\bf 7. Concluding remarks} \\

Of course, in such a short survey it is impossible to touch all the interesting
related topics. Let us make just a few remarks. \\

\noin {\bf 7.1.} Due to a lemma of T. Fujita [Fu 2, (2.4)] any smooth acyclic
algebraic surface is affine. In general, this does not hold in higher
dimensions. Indeed, J. Winkelmann [Wi] constructed a free regular
$\C_+$--action on $\C^5$ with the quotient $\C^5 // \C_+ = Q \setminus Z$,
where $Q$ is a smooth affine quadric of complex dimension four and $Z \subset
Q$ is a smooth codimension two subvariety. This quotient is diffeomorphic to
$\R^8$, but it is not Stein. \\

\noin {\bf 7.2.} By the Gurjar--Shastri Theorem [GuSha] any smooth acyclic
surface is rational. All the exotic $\C^n$--s in the present paper are rational
as well. For instance, Russell's threefolds are rational being
$\C^*$--varieties with a rational quotient. The general problem (posed by Van
de Ven [VdV]) {\it whether a smooth contractible quasiprojective variety is
rational}, is still open (for this and the related Hirzebruch problem on
compactifications of $\C^n$ see e.g. [MS, Fur]). \\

\noin {\bf 7.3.} At last, we remind the well known Abhyankar--Sathaye problem
on equivalence of embeddings $\C^k \hookrightarrow \C^n$, where $k < n \le
2k+1$. Note that for $n \ge 2k+2$ all such embeddings are equivalent
(Jalonek--Kaliman--Nori--Srinivas; see e.g. [Ka 4, Sr] and references therein),
while this is unknown already for the embeddings $\C \hookrightarrow  \C^3$ and
$\C^2 \hookrightarrow \C^3$. \\

\newpage

\begin{center} {\bf BIBLIOGRAPHY}\\
\end{center} \\
{\footnotesize

\noin [AM] S.S. Abhyankar, T.T. Moh, {\em Embedding of the line in the plane},
J. Reine Angew. Math. {\bf 276} (1975), 148--166.\\
\noin [BaDi] G. Barthel, A. Dimca, {\em On complex projective hypersurfaces
which are homology $\P_n$'s}, preprint MPI (1989) \\
\noin [ChoDi] A.D.R. Choudary, A. Dimca, {\em Complex hypersurfaces
diffeomorphic to affine spaces}, Kodai Math. J. {\bf 17} (1994), 171--178.\\
\noin [Di] A. Dimca, {\em Hypersurfaces in ${\bf C}^{2n}$ diffeomorphic to
${\bf R}^{4n - 2} \,(n \geq 2)$}, Max-Plank Institute, preprint, 1990.\\
\noin [FlZa] H. Flenner, M. Zaidenberg, {\sl Q-acyclic surfaces and their
deformations}, Proc. Conf. "Classification of Algebraic Varieties", Mai 22--30,
1992, Univ. of l'Aquila, L'Aquila, Italy /Livorni ed. Contempor. Mathem. {\bf
162}, Providence, RI, 1994, 143--208.\\
\noin [Fu 1] T. Fujita, {\em On Zariski problem}, Proc. Japan Acad. Ser. A
Math. Sci. {\bf 55} (1979), 106--110.\\
\noin [Fu 2] T. Fujita, {\em On the topology of non complete algebraic
surfaces}, J. Fac. Sci. Univ. Tokyo, Sect.IA, {\bf 29} (1982), 503--566.\\
\noin [Fur] M. Furushima, {\em The complete classification of compactifications
of $\C^3$ which are projective manifolds with second Betti number equal to
one}, Math. Ann. {\bf 297} (1993), 627--662.\\
\noin [GuMi] R.V. Gurjar, M. Miyanishi, {\sl Affine surfaces with $\k \le 1$},
Algebraic Geometry and Commutative Algebra, in honor of M. Nagata (1987),
99--124.\\
\noin [GuSha] R.V. Gurjar, A.R. Shastri, {\em On rationality of complex
homology 2--cells: I, II}, J. Math. Soc. Japan {\bf 41} (1989), 37--56,
175--212. \\
\noin [Ii 1] S. Iitaka, {\em On logarithmic Kodaira dimensions of algebraic
varieties}, Complex Analysis and Algebraic geometry, Iwanami, Tokyo, 1977,
175--189.\\
\noin [Ii 2] S. Iitaka, {\em Some applications of logarithmic Kodaira
dimensions}, Algebraic Geometry (Proc. Intern. Symp. Kyoto 1977), Kinokuniya,
Tokyo, 1978, 185--206.\\
\noin [IiFu] S. Iitaka, T. Fujita, {\em Cancellation theorem for algebraic
varieties}, J. Fac. Sci. Univ. Tokyo, Sect.IA, {\bf 24} (1977), 123--127.\\
\noin [Ka 1] S. Kaliman, {\em Smooth contractible hypersurfaces in ${\bf
C}^{n}$ and exotic algebraic structures on ${\bf C}^{3}$}, Math. Zeitschrift
{\bf
214} (1993), 499--510. \\
\noin [Ka 2] S. Kaliman, {\em Exotic analytic structures and
Eisenman intrinsic measures}, Israel Math. J. {\bf 88}(1994), 411--423.\\
\noin [Ka 3] S. Kaliman, {\em Exotic structures on $\C^n$ and $\C^*$-action
on $\C^3$}, Proc. on "Complex Analysis and Geometry", Lect. Notes in
Pure and Appl. Math., Marcel Dekker Inc. (to appear).\\
\noin [Ka 4] S. Kaliman, {\em Isotopic embeddings of affine algebraic varieties
into $\C^n$}, Contempor. Mathem. {\bf 137} (1992), 291--295. \\
\noin [KaML 1] S. Kaliman, L. Makar-Limanov,
{\em On some family of contractible hypersurfaces in ${\bf C}^4$ },
Seminair d'alg\`ebre. Journ\'ees Singuli\`eres et Jacobi\'ennes
26--28 mai 1993, Pr\'epublication de l'Institut Fourier, Grenoble,
(1994), 57--75. \\
\noin [KaML 2] S. Kaliman, L. Makar-Limanov, {\em On
Russell's contractible threefolds}, preprint, 1995, 1--22.\\
\noin [KaML 3] S. Kaliman, L. Makar-Limanov, {\em Affine algebraic manifolds
without dominant morphisms from Euclidean spaces}, preprint, 1995, 1--9 (to
appear in Rocky Mountain J.). \\
\noin [Kaw] Y. Kawamata, {\em Addition formula of logarithmic Kodaira dimension
for morphisms of relative dimension one}, Proc. Intern. Sympos. Algebr. Geom.,
Kyoto, 1977, Kinokuniya, Tokyo, 1978, 207--217. \\
\noin [KoRu] M. Koras, P. Russell, {\em On linearizing "good" ${\bf
C}^*$-action on ${\bf C}^3$}, Canadian Math. Society Conference Proceedings,
{\bf 10} (1989), 93--102.\\
\noin [Ko] M. Koras, {\em A characterization of ${\bf A}^2 / \Z_a$}, Compositio
Math. {\bf 87} (1993), 241--267.\\
\noin [Kr 1] H. Kraft, {\em Algebraic automorphisms of affine space},
Topological Methods in Algebraic Transformation Groups, Birkh\"auser, Boston
e.a., 1989, 81--105.\\
\noin [Kr 2] H. Kraft, {\em $\C^*$--actions on affine space}, Operator Algebras
etc., Progress in Mathem. {\bf 92}, 1990, Birkh\"auser, Boston e.a.,
561--579.\\
\noin [Lib] A. Libgober. {\em A geometric procedure for killing the middle
dimensional homology groups of algebraic hypersurfaces.} Proc. Amer. Math. Soc.
{\bf 63} (1977), 198--202.\\
\noin [LiZa] V. Lin, M. Zaidenberg, {\sl An irreducible simply connected curve
in ${\bf C}^{2}$ is equivalent to a quasihomogeneous curve}, Soviet Math.
Dokl., {\bf 28} (1983), 200-204. \\
\noin [ML] L. Makar-Limanov, {\em On the hypersurface $x+x^2y+z^2+t^3=0$
in $\C^4$}, preprint, 1994, 1--10. \\
\noin [MiSu] M. Miyanishi, T. Sugie, {\em Affine surfaces containing
cylinderlike open sets}, J. Math. Kyoto Univ., {\bf 20} (1980),
11--42.\\
\noin [Mi] M. Miyanishi, {\em Algebraic characterization of the affine
3--space}, Proc. Algebraic Geom. Seminar, Singapore, World Scientific, 1987,
53--67.\\
\noin [MS] S. M\"uller--Stach, {\em Projective compactifications of complex
affine varieties}, London Math. Soc. Lect. Notes Ser. {\bf 179} (1991),
277--283.\\
\noin [Ne] A. N\'emethi, {\em Global Sebastiani--Thom theorem for polynomial
maps}, J. Math. Soc. Japan, {\bf 43} (1991), 213--218.\\
\noin [PtD 1] T. Petrie, T. tom Dieck, {\em Contractible affine surfaces
of Kodaira dimension one}, Japan J. Math. {\bf 16}(1990), 147--169.\\
\noin [PtD 2] T. Petrie, T. tom Dieck, {\em The Abhaynkar--Moh problem in
dimension 3}, Lect. Notes Math. 1375, 1989, 48--59.\\
\noin [PtD 3] T. Petrie, T. tom Dieck, {\em Homology planes. An announcement
and survey}, Topological Methods in Algebraic Transformation Groups, Progress
in Mathem., {\bf 80}, Birkhauser, Boston, 1989, 27--48.\\
\noin [Pe] T. Petrie, {\em Topology, representations and equivariant algebraic
geometry}, Contemporary Math. {\bf 158}, 1994, 203--215.\\
\noin [Ram] C.P. Ramanujam, {\sl A topological characterization of the affine
plane as an algebraic variety}, Ann. Math., 94 (1971), 69-88.\\
\noin [Ru 1] P. Russell, {\em On a class of ${\bf C}^3$-like threefolds},
Preliminary Report, 1992.\\
\noin [Ru 2] P. Russell, {\em On affine-ruled rational surfaces},
Math. Ann., {\bf 255} (1981), 287--302.\\
\noin [Sa] F. Sakai, {\em Kodaira dimension of complement of divisor}, Complex
Analysis and Algebraic geometry, Iwanami, Tokyo, 1977, 239--257.\\
\noin [Sr] V. Srinivas, {\em On the embedding dimension of an affine variety},
Math. Ann.  {\bf 289} (1991), 125-132. \\
\noin [Su] M. Suzuki, {\em Propi\'{e}tes topologiques des polyn\^{o}mes de
deux variables complexes, et automorphismes alg\'{e}brigue de l'espace
${\bf C}^{2}$}, J. Math. Soc. Japan, {\bf 26} (1974), 241-257. \\
\noin [tD 1] T. tom Dieck, {\em Hyperbolic modifications and acyclic
affine foliations}, preprint, Mathematisches Institut, G\"ottingen, H.27
(1992), 1--19.\\
\noin [tD 2] T. tom Dieck, {\em Ramified coverings of acyclic varieties},
preprint, Mathematisches Institut, G\"ottingen, H.26
(1992), 1--20.\\
\noin [VdV] A. Van de Ven, {\em Analytic compactifications of complex homology
cells}, Math. Ann. {\bf 147} (1962), 189--204.\\
\noin [Wi] J. Winkelmann, {\em On free holomorphic $\C^*$--actions on $\C^n$
and homogeneous Stein manifolds}, Math. Ann. {\bf 286} (1990), 593--612.\\
\noin [Za 1] M. Zaidenberg, {\sl Ramanujam surfaces and exotic algebraic
structures on ${\bf C}^n$}, Soviet Math. Doklady {\bf 42} (1991), 636--640.\\
\noin [Za 2] M. Zaidenberg, {\sl An analytic cancellation theorem and exotic
algebraic structures on ${\bf C}^n$, $n \ge 3$},  Ast\'erisque {\bf 217}
(1993), 251--282. \\
\noin [Za 3] M. Zaidenberg, {\sl Isotrivial families of curves on affine
surfaces and characterization of the affine plane}, Math. USSR Izvestiya {\bf
30} (1988), 503-531. Addendum: ibid, {\bf 38} (1992), 435--437.\\
\noin [Za 4] M. Zaidenberg, {\sl On Ramanujam surfaces, $\C^{**}$-families
and exotic algebraic structures on $\C^n$, $n \ge 3$},
Trudy Moscow Math. Soc. {\bf 55} (1994), 3--72 (Russian; English transl.
to appear). \\

\noin  Mikhail Zaidenberg

\noin Universit\'{e} Grenoble I

\noin Laboratoire de Math\'ematiques

\noin associ\'e au CNRS

\noin BP 74

\noin 38402 St. Martin d'H\`{e}res--c\'edex

\noin France

\noin e-mail: ZAIDENBE@FOURIER.GRENET.FR}

\end{document}